file
the
pict.uu '' .
\newcommand{\beq}{\begin{equation}}
\newcommand{\eq}{\end{equation}}
\newcommand{\barr}{\begin{array}}
\newcommand{\earr}{\end{array}}
\newcommand{\beqa}{\begin{eqnarray}}
\newcommand{\eqa}{\end{eqnarray}}
\newcommand{\baqa}{\begin{eqnarray*}}
\newcommand{\aqa}{\end{eqnarray*}}
\newcommand{\defi}{\newcommand}
\newcommand{\hook}{\hookrightarrow}

\defi{\noi}{\noindent}
\defi{\NN}{I\!\!N}
\defi{\CC}{C\!\!\!\!I}
\defi{\QQ}{Q\!\!\!\!I}
\defi{\ZZ}{Z\!\!\!Z}
\defi{\CN}{{\bf C}$\!\!\!$l\quad}
\defi{\RN}{{\bf R}^+$\!\!\!\!\!$l}

\newlength{\partcolwidth}\newlength{\restcolwidth}

\defi{\A}{{\cal A}}
\defi{\AS}[4]{  {\hat A_q^{#1#2}}{}_{#3#4}  }
\defi{\sis}[4]{  {\hat S_q^{#1#2}}{}_{#3#4}  }
\defi{\SSS}[4]{  {\hat S_{\sqrt{q}}^{#1#2}}{}_{#3#4}  }
\defi{\TS}[4]{  {\hat T_q^{#1#2}}{}_{#3#4}  }
\defi{\C}{{\bf C}}
\defi{\E}{{\bf 1}}
\defi{\M}[2]{  {M^{#1}}{}_{#2}  }
\defi{\meta}[2]{  {\eta^{#1}}{}_{#2}  }
\defi{\m}[2]{  {m^{#1}}{}_{#2}  }
\defi{\CG}[3]{ c^{#1}_{#2#3}}
\defi{\MM}{{M}}
\defi{\POP}{{\cal P}}
\defi{\RA}{\hat{R}}
\defi{\ra}{\hat{r}}
\defi{\RMA}[4]{  {\hat R_q^{#1#2}}{}_{#3#4}  }
\defi{\RMAS}[4]{  {\hat R_{\sqrt{q}}^{#1#2}}{}_{#3#4}  }
\defi{\rma}[4]{  {\hat r^{#1#2}}{}_{#3#4}  }
\defi{\Ri}[4]{  {\hat R^{-1}_q{}^{#1#2}}{}_{#3#4}  }
\defi{\RR}[4]{  R_q^{#1#2}{}_{#3#4}  }
\defi{\RRi}[4]{  R^{-1}_q{}^{#1#2}{}_{#3#4}  }
\defi{\war}{\bar{w}}
\defi{\wi}{w}
\defi{\SU}{SU_{\sqrt{q}}(2)}
\defi{\SUX}{SU^{ex}_{\sqrt{q}}(2)}
\defi{\SOd}{SO_q(3)}
\defi{\SOz}{SO_q(2,1)}
\defi{\SL}{SL_q(2,R)}
\defi{\SLX}{SL^{ex}_q(2,R)}
\defi{\nomb}{\nonumber}
\defi{\al}{\alpha}
\defi{\be}{\beta}
\defi{\bx}{\bar{x}}
\defi{\bv}{\bar{v}}
\defi{\frq}{\frac{1}{Q}}
\defi{\ga}{\gamma}
\defi{\de}{\delta}
\defi{\eps}{\epsilon}
\defi{\veps}{\varepsilon}
\defi{\Th}{\Theta}
\defi{\Thi}{\Theta_i}
\defi{\vth}{\vartheta}
\defi{\ka}{\kappa}
\defi{\la}{\lambda}
\defi{\vrh}{\varrho}
\defi{\roi}{\rho_i}
\defi{\si}{\sigma}
\defi{\vph}{\varphi}
\defi{\om}{\omega}
\defi{\bom}{\bar{\om}}

\documentstyle[a4,11pt]{article}

\begin{document}

\makebox[13.5cm][r]{LMU--TPW 95--6}
\vspace{3cm}
\Large

\begin{center}
{\bf $ISO_q(3)$ and $ISO_q(2,1)$}

\vspace{2cm}

{\large {\it Sebastian Sachse, Ralf Weixler\\ \it Sektion Physik der
Universit\"at M\"unchen
\\
\it Lehrstuhl Professor Wess, Theresienstra\ss e 37\\
\it D-81245 M\"unchen , Germany }}

\vspace{2cm}
{\large \today}
\end{center}
\vspace{2cm}

%\title{$ISO_q(3)$ and $ISO_q(2,1)$}
%\author{Sebastian Sachse, Ralf Weixler}
%\maketitle
\normalsize
\begin{abstract}
We prove the embedding of $ ISO_q(3) \hook ISU^{ex}_{\sqrt{q}}(2) $
and $ ISO_q(2,1) \hook ISL^{ex}_q(2,R)$ as $^*$-algebras and give a Hilbert
space
representation of $I\SUX.$
\end{abstract}
\section{ Introduction}

The inhomogenized extensions of a large list of standard quantizied
Lie groups \cite{frt} have been given in \cite{weix1,weix2,kob,drab,fiore,maj}.
They
form quantizied versions of the classical inhomogeneous groups.
For a real deformation parameter $q$
the representation theory of the homogeneous parts (e.g. corepresentations
of the function algebra) is basically the same as for the classical
groups, whereas for $q$ root of unity it is completely different.
The representation theory of the noncommutative function algebra however
differs for any $q \neq 1$ from the classical situation. Its relevance
stems from the question whether a deformation exists on the $C^*$-algebra
level. \cite{woro}

In part 2 we recall the properties of inhomogeneous quantum groups.
In the 3rd. part we examine the algebraic embedding of the
$ISO_q(3)$  into $ISU^{ex}_{\sqrt{q}}(2)$ and $ISO_q(2,1)$ into $I\SLX$. Here
the
``extended inhomogeneous'' quantum algebra  $IG^{ex}$ designates inhomogeneous
quantum algebras containing two sets of coordinate functions.

In the last chapter we examine the representation theory of the  $I\SUX$.

\section{ The Hopf algebra structure of inhomogeneous quantum groups }

Quantum groups may be considered to be deformations of the function algebra
over the corresponding Lie groups.
The deformation is given by a parameter $q \epsilon$ \CN
which has to be further restricted in order to get special cases of
deformations.
Quantum groups exhibit a Hopf algebra structure.
The noncommutative algebra structure is controlled by an $\RA$--matrix
fulfilling the
Quantum Yang--Baxter equation $\RA_{12} \RA_{23} \RA_{12} = \RA_{23} \RA_{12}
\RA_{23}$.
In this paper we refer to $\RA$--matrices in their standard form given in
\cite{frt}.
They are defined by their projector decomposition making use of the
antisymmetrizer $\AS ijkl $, the symmetrizer $\sis ijkl $ and the trace
projector $\TS ijkl\propto C^{ij} C_{kl}$ with the metric $C_{ij}$,
existing for the $q$--orthogonal groups only.
\beqa
\RMA ijkl = \left\{ \barr {ll}
 q\sis ijkl - q^{-1}\AS ijkl & \qquad(\mbox{ for } SL_q(N)) \\
           q\sis ijkl - q^{-1}\AS ijkl + q^{1-N} \TS ijkl & \qquad(\mbox{ for }
SO_q(N))

\label{1.1}
\earr
\right.
 \eqa

The algebra relations for the generators $\M ij$ of the unital
\CN-algebra $\A$ are:

\beq
\RMA ij{j'}{i'} \M {j'}{j''} \M {i'}{i''} =
   \M i{i'} \M j{j'} \RMA {i'}{j'}{j''}{i''}
   \label{1.2}
\eq
and
\beq
 \left\{ \barr  {ll}
     \det M = \frac{(-1)^{N-1}}{[N]!}\epsilon^{k_1...k_N}
    \M {l_1}{k_1} ... \M {l_N}{k_N}\epsilon_{l_1...l_N} = \E & (SL_q(N)) \\
    C_{ij} \M i{i'} \M j{j'} = C_{i'j'}\E & (SO_q(N))
\earr \right.    \label{1.3}
\eq

For the unimodularity condition we use the $q$--antisymmetric tensors
$\epsilon_q$ defined in \cite{cssw}.

In order to obtain inhomogeneous quantum groups the set of generators
has to be enlarged not only by the coordinate functions $x_i$ but by an
invertible scaling operator $\war$ as well. Its existence is required  by
consistency of the comultiplication. The additional algebra relations of the
extended Hopf algebra
$\A^I$ are:

\parbox[t]{6cm}{\baqa
  && (i) \quad
      x^i \M jk = \gamma \RMA ijlm \M lk x^m \\
  && (ii) \quad
      \war \wi = \E \\
 && (iii) \quad
      \war \M ij = \M ij\war
\aqa}
\hfill
\parbox[t]{6cm}{\baqa
&& (iv) \quad    \wi\M ij = \M ij\wi   \\
  && (v) \quad   \war x^i = \frac{q}{\gamma} x^i\war  \\
   && (vi) \quad  \wi x^i = \frac{\gamma}{q} x^i\wi
\aqa}
\hfill
\begin{minipage}[t]{1cm}
\vspace{1cm}
\beq  \label{1.4} \eq
\end{minipage}
\nopagebreak
with
\beq
     \gamma = \left\{\barr{ll} q^{-1/N} & \qquad (\mbox{ for } SL_q(N)) \\
                     1              & \qquad (\mbox{ for } SO_q(N)) \earr
\right.
\eq
The comultiplication $\Phi: \A^I \to \A^I \otimes \A^I $, counit
$e: \A^I \to$ \CN  and the antipode $\kappa : \A^I \to \A^I $
are very easily given in matrix notation.

With
\beq
M^I = \left( \barr {cc}
\war M & {\bf x} \\ 0 & 1 \earr \right)
\eq
we get
\beqa
 \
\Phi(\MM^I) = \MM^I \dot{\otimes} \MM^I
\eqa
and
\beqa
  \hspace{1.5cm}  e(\MM^I) =
\left( \barr {cc}
E & {\bf 0} \\ 0 & 1 \earr \right)
\eqa
with the unity matrix $E$.

The antipode is given as
\beq
\kappa(\MM^I) = \left( \barr {cc}
\kappa(\MM)\war & -\kappa(\MM)\war\vec x \\ 0  & \E  \earr \right).
\eq

\subsection{ Complex conjugation}

The $*$-operations on $\A $ are defined quite differently in the cases  $ q
\epsilon \mbox{\bf R}^+ $ and $ |q|=1 $.
a) $ q \epsilon \mbox{ \bf R}^+$

With the unitarity condition $ ( \M ij )^* = \kappa (\M ji )$, the
quantum group $SL_q(N)$ becomes a $SU_q(N)$ quantum group. The same
*-structure holds for the orthogonal quantum groups $SO_q(N)$ .

b) $|q| = 1$

For such $q$ the $R_q$ matrix has the property $R_q^* = R_q^{-1}$.
With the reality
condition $ ( \M ij )^* = \M ij $ one finds the real representation of
the quantum group $SL_q(N) $ called $SL_q(N,R)$ and the orthogonal
quantum groups in this case have a metric which is indefinite, i.e. for
$N$ even we get $SO_q(n,n)$ and for $N$ odd $SO_q(n,n+1)$, with $N = 2n$ or $N
=
2n+1$ respectivly.

The complex conjugation for the inhomogeneous extensions of these
function algebras have to be treated separatly as well.

a) In the case $ q \epsilon \mbox{ \bf R}^+$ we have to enlarge the
generating set of the $^*$-Hopf algebra $\A^I$ by the conjugate coordinate
functions $\bar{x}_i$.
The additional algebra relations are:

\parbox[t]{6cm}{\baqa
&& (i)\quad
      \M ls \bar{x}_j = \frac{1}{\gamma}\RMA alij \bar{x}_a\M is  \\
&&  (ii)\quad
      \wi \bar{x}_i = \frac{\gamma}{q} \bar{x}_i \wi
\aqa}
\hfill
\parbox[t]{6cm}{\baqa
&& (iii)\quad
      \war \bar{x}_i = \frac{q}{\gamma} \bar{x}_i \war  \\
&& (iv)\quad
      x^i \bar{x}_j =  \frac{1}{q} \bar{x}_a x^b \RMA aibj.
\aqa}
\hfill
\begin{minipage}[t]{1cm}
\vspace{1cm}
\beq  \eq
\end{minipage}

The comultiplication of $\bar{x}_i$ follows from being an
*-homomorphism and the antipode from the fact that $\kappa \circ *
\circ \kappa \circ * = $ id.

b) When $|q| = 1$ the coordinate functions may be chosen to be real $(x_i^* =
x_i)$, since applying the * to (\ref{1.4}(i)) and taking into
account that
$\RR ijkl ^* = \RRi ijkl $,
we get:
\beq
       \M jk x^i = \gamma^{-1} \RRi jilm x^m \M lk  \nomb
\eq
or
\beq
    \gamma \RMA rsji \M jk x^i = \delta^r_m \delta^s_l x^m \M lk .
\eq

\section{ Algebraic embedding}

The $\RA -$matrix of $\SU$ is decomposed using the $q$-antisymmetric
$\veps$-tensor:
\beq
\RMAS \mu \nu \rho \si = \sqrt{q} \delta_\rho^\mu \delta_\nu^\si +
\veps^{\mu \nu} \veps_{\rho \si},
\eq
with

\[\eps_{\mu\nu} := \left( \barr{cc} 0&q^{-1/4}\\-q^{1/4}&0 \earr
\right)_{\mu\nu} = -\eps^{\mu\nu}
 \]

For the homogeneous parts the embedding $ SO_q(3) \hook SU_{\sqrt{q}}$ is well
known \cite{schl}. With the q-deformed Clebsch-Gordan-coefficients $\CG i
\mu\nu$ of the
product decomposition of the $\SU$\cite{weich}
%\beqa
%\CG i 11 & = & (1,0,0) \nomb \\
%\CG i 12 & = & \frac{\sqrt{q}}{\sqrt{1+q}} (0,1,0)  \\
%\CG i 21 & = & \frac{1}{\sqrt{1+q}} (0,1,0) \nomb \\
%\CG i 22 & = & (0,0,1) \nomb
%\eqa
\beq
\CG 1 \mu\nu := \left( \barr{cc} 1&0\\0&0 \earr \right)_{\mu\nu},\qquad
\CG 2 \mu\nu := \frac{1}{\sqrt{1+q}}\left( \barr{cc} 0&\sqrt{q}\\1&0 \earr
 \right)_{\mu\nu},\qquad
\CG 3 \mu\nu := \left( \barr{cc} 0&0\\0&1 \earr \right)_{\mu\nu}\quad
\eq
and $\CG i \mu\nu = c_i^{\mu\nu}$ the matrix elements $\M ij$ of $\SOd$ are
given in terms of $\SU$ elements $\m \mu\nu$ by
\beq
\M ij := \CG i \mu\nu \m \mu\rho \m \nu\sigma c_j^{\rho\sigma}.
\eq
This is  verified by checking the $\RA$-matrix of the $\SOd$ group to be :
\beq
\frac{1}{q} \meta i{i'} \meta j{j'} \CG i \mu\rho \CG j \sigma\lambda
\rma \mu\nu {\mu'}{\nu'} \rma \rho\sigma \nu\tau \rma \tau\lambda
{\tau'}{\lambda'} \rma {\nu'}{\tau'} {\rho'}{\sigma'}
c_k^{\mu'\rho'} c_l^{\sigma'\lambda'}, \label{2.4}
\eq
where $\ra$ denotes the $\RA$-matrix of the $\SU$ group and
$\eta = \mbox{diag}(1,i,1)$, $i$ being the imaginary unit. It just produces a
base change to more convenient
coordinates.
Note as well the decomposition of the symmetric
projector $\SSS \mu \nu \rho \si = c_i^{\mu \nu} \CG i \rho \si.$
\noi

In order to clear out nasty indices we want to
make use of a graphical notation, which has been given in \cite{cssw}.
With
\vspace{1.5cm}

\parbox{13cm}{
\beq\hspace{0cm}\CG i \mu\nu := \hspace{1.2cm}\eq}
\vspace{-2cm}

\hspace{7cm}
\begin{minipage}{5cm}
\input{epsf}
\epsffile{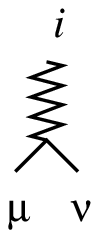}
\end{minipage}
and the equality
\vspace{.5cm}

\hspace{2cm}
\begin{minipage}{5cm}
\input{epsf}
\epsffile{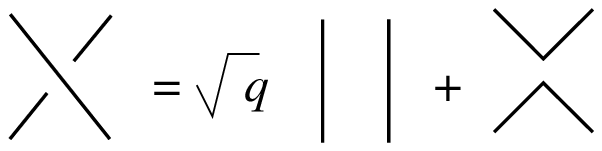}
\end{minipage}
\vspace{-1cm}
\beq \hspace{.5cm}    \eq
\vspace{.5cm}
\noi

we can disentangle the matrix:
\vspace{.5cm}

\hspace{0cm}
\begin{minipage}{5cm}
\input{epsf}
\epsffile{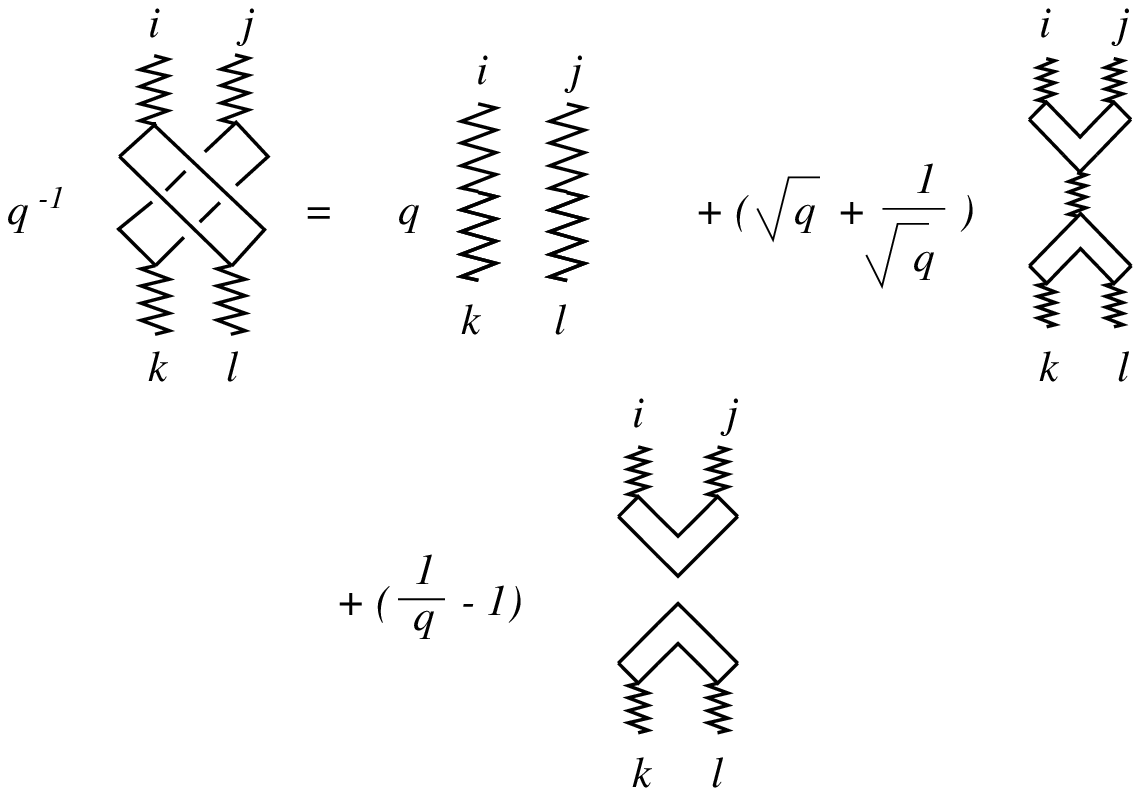}
\end{minipage}
\nopagebreak
\vspace{-3cm}
\beq \hspace{.5cm}    \eq
\vspace{2cm}

This is the $\RA$-matrix of $\SOd$ since obviously
\vspace{1cm}

\parbox{13cm}{
\beq\hspace{-1cm}C^{ij} = \meta i{i'} \meta j{j'}\hspace{1.2cm}\eq}
\vspace{-1.5cm}

\hspace{7cm}
\parbox{4cm}{
\input{epsf}
\epsffile{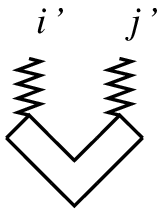}}
and
\vspace{2.5cm}

\parbox{12cm}{\[\hspace{2.5cm}\AS ij kl = \meta i{i'} \meta j{j'}\hspace{3cm}
\meta {k'}k
\meta {l'}l\hspace{1.2cm} \]}
\parbox{1cm}{\beq \eq}
\vspace{-3.5cm}

\hspace{6.5cm}
\parbox{10cm}
{\input{epsf}
\epsffile{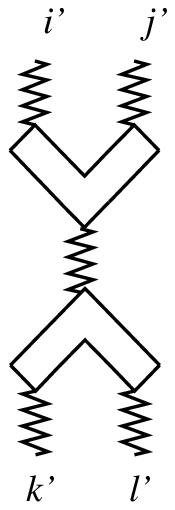} }
\vspace{.5cm}

Of course the same construction holds for the $\SOz$ group. Then $\eta$ has
to be chosen as identity matrix.

b) Since we know the $q$-antisymmetrizer we are able to find the $\SOd$
covariant
quantum plane in terms of spinor variables. To obtain sufficiently many
degrees of freedom we have to take at least two copies of $q$-spinors
 $x$ and $y$ having the same commutation relations with $m$.
This provides an extended inhomogeneous algebra called $I\SUX$.
We want to mention that the extended algebra does not have a correct coalgebra
structure.
This is not important for the algebraic embedding.
The 3-dimensional quantum space has the form:
\vspace{1cm}

\beq \hspace{-5cm} z^i = \meta i {j} \eq
\vspace{-2cm}

\hspace{5cm}
\begin{minipage}{8cm}
\input{epsf}
\epsffile{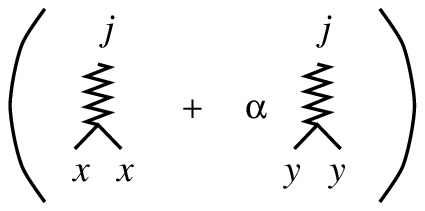}
\end{minipage}
\vspace{.7cm}

We fix the $x,y$- relations such that the element $\varepsilon_{\nu \mu} \,
x^{\nu} y^{\mu} $ commutes with the
coordinate functions and
get
\vspace{.5cm}

\beq \hspace{.5cm}  =  \eq
\nopagebreak
\vspace{-1.7cm}

\hspace{5cm}
\begin{minipage}{8cm}
\input{epsf}
\epsffile{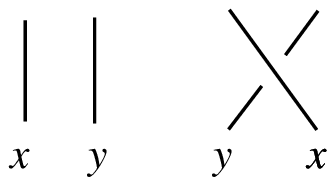}
\end{minipage}
\noi

We still have to prove that $ \AS ijkl z^i z^j $vanishes.
This follows from the equation
\vspace{1cm}

\beq\hspace{4cm} = 0,\eq
\nopagebreak
\vspace{-2cm}

\hspace{3cm}
\begin{minipage}{4cm}
\input{epsf}
\epsffile{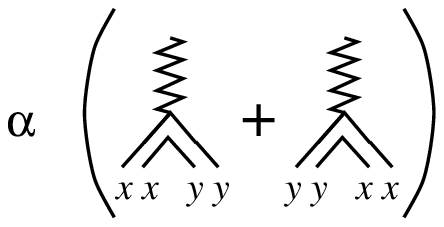}
\end{minipage}
\noi

or
\vspace{1.5cm}

\beq  \eq
\vspace{-2cm}

\hspace{2cm}
\begin{minipage}{8cm}
\input{epsf}
\epsffile{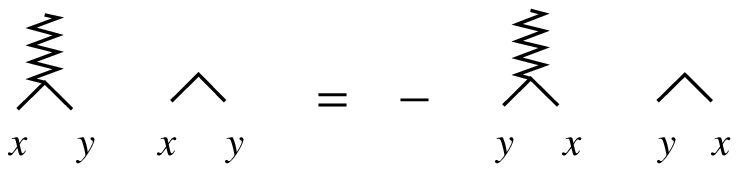}
\end{minipage}
\noi

Next we observe that

\hspace{1.5cm}
\parbox{10cm}
{\input{epsf}
\epsffile{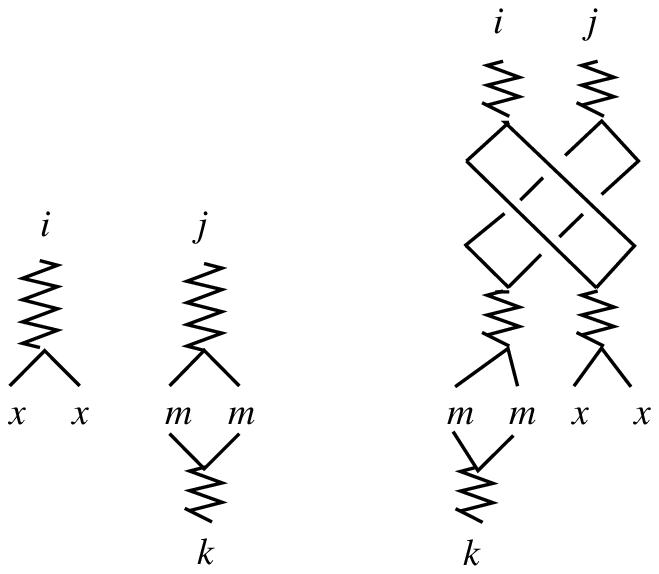} }
%\vspace{-2.5cm}
\vspace{-2.8cm}

\beq
\hspace{1.5cm}
%x\hspace{.5cm} x \hspace{.8cm} m\hspace{.5cm} m  \hspace{.5cm} = q^{-1}
\hspace{.4cm}m\hspace{.5cm} m \hspace{.4cm}x
%\quad x
\eq
\vspace{2cm}

Now we have to take a look to the * - operation.
The behaviour is for $q$ complex quite different from $q$ real,
since $(\CG i \mu\nu )^* = \CG i \nu\mu $ in the first and
$(\CG i \mu\nu )^* = \CG i \mu\nu $ in the second case.
{}From this we see immediately that the algebraic embedding in the case of
$I\SOz $ respects the * - operation,
whereas for the $I\SOd$ quantum group we have to examine the * - structure
more closely.
At first for the homogeneous part we have:
\beq
\kappa (\CG i \mu\nu \m \mu\rho \m \nu\sigma c_j^{\rho\sigma})
= \CG i \mu \nu (\m \si\nu )^* (\m \rho\mu)^* c_j^{\rho\si }
= \left(\CG j \rho\si \m \rho\mu \m \si\nu c_i^{\mu\nu} \right)^*.
\eq
To use the graphical technique for the translations
we have to introduce the ``transposed $c$--matrix ''$c^T$:
\vspace{1cm}

\hspace{4cm}
\begin{minipage}{8cm}
\input{epsf}
\epsffile{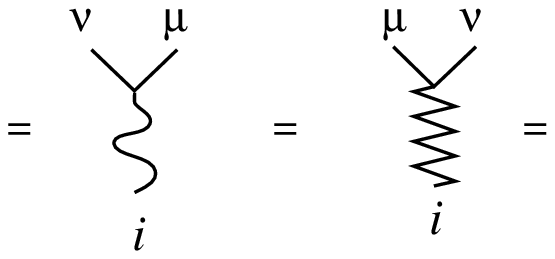}
\end{minipage}
\nopagebreak
\vspace{-2.2cm}

\beq
 (c^T)^{\nu\mu}_i \hspace{7cm} c_i^{\mu\nu}
\eq
\vspace{1.5cm}
\noi

Of course $(z^i)^* =: \bar{z}_i = (c_i^T)^{\mu\nu} \bar{x}_\mu \bar{x}_\nu $
Calculating now the mixed commutation relations one has:

\hspace{0cm}
\begin{minipage}{5cm}
\input{epsf}
\epsffile{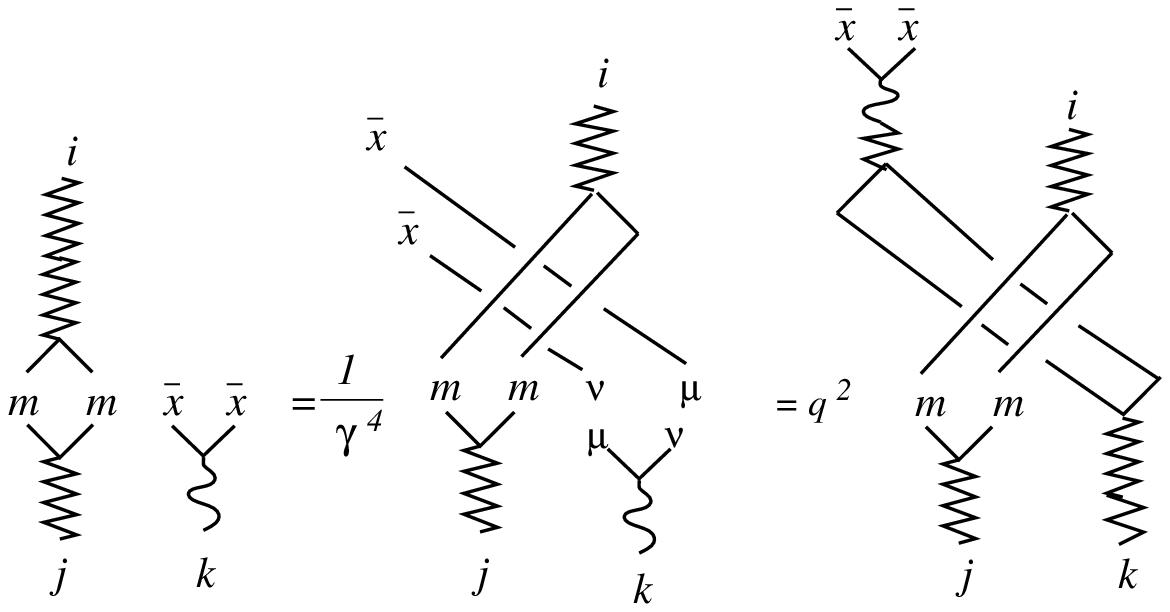}
\end{minipage}
\vspace{-3cm}
\beq \hspace{.5cm}  \label{29}  \eq
\vspace{2cm}
\noi

To reproduce the correct $z-\bar{z}$ relations
which use diagramms similar to that of (\ref{29}) we have to redefine
$z$ as $z \rightarrow z^i = x^\mu x^\nu \CG i \mu\nu \bar{\om}^n $. Now the
calculation is similar to the one above and we find $n= 2/3$.
This relation finishes the prove of the algebraic embedding.

\noi

Remark:

a) The antipode of the coordinates $z$ just differs by a minus sign, when
expressed in terms of $\ka(x)$.

b) The Coproduct of the coordinates $z$ can't be embedded by obvious reasons.

\section{Hilbert space representation for $I\SOd$}

The problem to find a representation for $I\SOd$ and $I\SOz$
is now reduced to that of representing $I\SUX$ and $I\SLX$.

Here we restrict ourselves to the representation of $I\SOd$. (We want to
mention
that $\SL$ does not exist on the Hilbert space level anyhow \cite{Schmud}.)
As well we confine the value of $q$ to $(0,1)$. The case $q>1$ is however
isomorphic. The relations of $\SU$ are ($\mu = \sqrt{q}$):
\beqa
\al \ga = \mu \ga \al, \qquad &\al \ga^* = \mu \ga^* \al&, \qquad \ga \ga^* =
\ga^* \ga, \nomb \\
\al^* \al + \ga^* \ga = 1, &&  \al \al^* + \mu^2 \ga^* \ga =1
\eqa
The second coordinate has the following commutation relations:
(Remember that in our convention the quantum plane coordinates are $x^{(\tau)},
\tau=1,2$.)
\beqa
x^{(2)} \al &= \mu^{-1} \al x^{(2)}, \qquad x^{(2)}  \ga =& \mu \ga x^{(2)}
\nomb \\
  x^{(2)} \al^* &= \mu \al^* x^{(2)}, \qquad x^{(2)}  \ga^* =& \mu^{-1} \ga^*
x^{(2)}  \label{4.1}
\eqa
Applying the antipode to (\ref{4.1}) gives
\beqa
\al^* \ka(x^{(2)}) &= \mu^{-1} \ka(x^{(2)}) \al^*, \qquad
 \ga \ka(x^{(2)}) =& \mu \ka(x^{(2)}) \ga \nomb \\
\al \ka(x^{(2)})  &= \,\mu \ka(x^{(2)}) \al, \qquad
\quad \ga^* \ka(x^{(2)}) =& \mu^{-1} \ka(x^{(2)}) \ga^*.
\eqa
Those relations hold for any $\SU$-covariant plane.
\noi

We have other relations for the functions on the
quantum planes $x^{(\tau)}$ and $y^{(\tau)}, \tau = 1,2$:

\parbox[t]{6cm}{\baqa
x^{(2)} y^{(2)} &=& \mu  y^{(2)} x^{(2)},   \\
x^{(\tau)} \ka(x^{(\de)}) &=& \mu^{-2} \ka(x^{(\de)}) x^{(\tau)},  \\
\ka(x^{(2)})  \overline{\ka(y^{(2)})} &=& \mu^2 \overline{\ka(y_2)} \ka(x_2),
\\
\overline{a_\nu} \ka(b^{(\tau)}) &=& \ka(b^{(\tau)}) \overline{a_\nu},
\aqa}
\hfill
\parbox[t]{6cm}{\baqa
\ka(x^{(2)}) \ka(y^{(2)}) &=& \mu^{-1} \ka(y^{(2)}) \ka(x^{(2)}),  \\
x^{(\tau)} \ka(y^{(\de)}) &=& \mu^{-1} \ka(y^{(\de)}) x^{(\tau)},   \\
 \\
a,b \,\eps\, \{x, y\},\hspace{.5cm}&&
\aqa}
\hfill
\begin{minipage}[t]{1cm}
\vspace{1cm}
\beq \eq
\end{minipage}
\vspace{-.5cm}
\[
y^{(\tau)} \ka(x^{(\de)}) = \mu^{-1} \ka(x^{(\de)}) y^\tau +
(\mu^{-2}-1) \ka(y^{(\de)}) x^{(\tau)},
\]
where $\de, \tau \,\eps\, \{1,2\}$. \\
These and the conjugated relations together with the
obvious relations for $\om$ contain the whole algebraic information of $I\SU$
since:
\beq a^{(1)} = (\mu \ga)^{(-1)} \left(\om \ka\left(
a^{(2)}\right) + \al a^{(2)} \right), \qquad a,b \,\eps\, \{x, y\}.
\eq
Remember that $\ga $ is an invertible element.

The algebra may still be simplyfied by a nonlinear transformation in the
functions of coordinates.
With $Q^4 = \mu^2 = q$ and $v^3 = \om $ we define:

\parbox[t]{6cm}{\baqa
\rho_1 &=& \bv^2 \ga^{-1} x^{(2)} \\
\rho_2 &=& \bv^2 \ga^{-1} y^{(2)}
\aqa}
\hfill
\parbox[t]{6cm}{\baqa
\Th_1 &=& q^{-1} v \ga^{-1} \ka(x^{(2)}) \\
\Th_2 &=& q^{-1} v \ga^{-1} \ka(y^{(2)})
\aqa}
\hfill
\begin{minipage}[t]{1cm}
\vspace{1cm}
\beq \eq
\end{minipage}
\noi

All algebraic relations with coordinate-functions are given by:

\parbox[t]{4cm}{\baqa
\roi \al &=& Q \al \roi \\
\roi \al^* &=& Q^{-1} \al^* \roi \\
\roi \ga &=& Q \ga \roi \\
\roi \ga^* &=&Q^{-1} \ga^* \roi
\aqa}
\hfill
\parbox[t]{4cm}{\baqa
\Thi \al &=& Q \al \Thi \\
\Thi \al^* &=& Q^{-1} \al^* \Thi \\
\Thi \ga &=& Q^{-1} \ga \Thi \\
\Thi \ga^* &=& Q \ga^* \Thi
\aqa}
\hfill
\parbox[t]{4cm}{\baqa
\Thi \Thi^* &=& \Thi^* \Thi \\
\roi \roi^* &=& Q^2 \roi^* \roi \\
\roi \Thi &=& Q^{-3} \Thi \roi \\
\roi^* \Thi &=& Q^{-1} \Thi \roi^*
\aqa}
\hfill
\begin{minipage}[t]{1cm}
\vspace{1cm}
\beq \eq
\end{minipage}
with $i=1,2$.

\parbox[t]{6cm}{\baqa
\Th_1 \Th_2 &=& Q^{2} \Th_2 \Th_1 \\
\rho_1 \rho_2 &=& Q^2 \rho_2 \rho_1 \\
\rho_1 \Th_2 &=& Q^{-1} \Th_2 \rho_1 \\
\Th_1 \rho_2 &=& (Q^2 - Q^{-2}) \Th_2 \rho_1 + Q \rho_2 \Th_1
\aqa}
\hfill
\parbox[t]{6cm}{\baqa
\rho_1 \rho^*_2 &=& Q^2 \rho_2^* \rho_1 \\
\Th_1 \rho_2^* &=& Q \rho_2^* \Th_1 \\
\Th_2 \rho_1^* &=& Q \rho_1^* \Th_2 \\
\Th_1 \Th_2^* &=& \Th_2^* \Th_1
\aqa}
\hfill
\begin{minipage}[t]{1cm}
\vspace{1cm}
\beq \eq
\end{minipage}
\noi

It is easy to find a maximal real subalgebra ${\cal D}$ of commuting elements.
In order to rely the representation of $I\SOd$ to that of $\SU$ we choose:
\beq
{\cal D}:= \{\al \al^*, \ga \ga^*, \rho_1 \rho_1^*, \Th_1 \Th_1^*, \Th_2
\Th_2^*, v\}.
\eq
These operators are used to label the eigenvectors of the Hilbertspace. They
are normalised so that the representation is given by:
\beqa
 \pi(\al) |n,m,k,r,s,v\rangle &=& \sqrt{1-Q^{4n}} |n-1,m,k,r,s,v\rangle \nomb
\\
 \pi(\al^*) |n,m,k,r,s,v\rangle &=& \sqrt{1-Q^{4(n+1)}} |n+1,m,k,r,s,v\rangle
\nomb \\
 \pi(\ga) |n,m,k,r,s,v\rangle &=& Q^{2n} |n,m-1,k,r,s,v\rangle \nomb \\
 \pi(\ga^*) |n,m,k,r,s,v\rangle &=& Q^{2n+2} |n,m+1,k,r,s,v\rangle \nomb \\
 \pi(v) |n,m,k,r,s,v\rangle &=& |n,m,k-1,r,s,v\rangle \nomb \\
 \pi(v^*) |n,m,k,r,s,v\rangle &=& |n,m,k+1,r,s,v\rangle \nomb \\
 \pi(\rho_1) |n,m,k,r,s,v\rangle &=& Q^{-(n+m+k)+r+2s-v}|n,m,k,r-1,s,v\rangle
\nomb \\
 \pi(\rho_1^*) |n,m,k,r,s,v\rangle &=&
Q^{-(n+m+k)+r+1+2s-v}|n,m,k,r+1,s,v\rangle \nomb \\
 \pi(\rho_2) |n,m,k,r,s,v\rangle &=&
Q^{-(n+m+k)+2r+3s-v}|n,m,k,r-1,s,v-1\rangle
\nomb \\
 \pi(\rho_2^*) |n,m,k,r,s,v\rangle &=&
Q^{-(n+m+k)+2r+1+3s-v}|n,m,k,r+1,s,v+1\rangle \nomb \\
 \pi(\Th_1) |n,m,k,r,s,v\rangle &=& Q^{-(n-m+k+r+v)}|n,m,k,r,s-1,v\rangle \nomb
\\
 \pi(\Th_1^*) |n,m,k,r,s,v\rangle &=&  Q^{-(n-m+k+r+v)}|n,m,k,r,s+1,v\rangle
\nomb \\
 \pi(\Th_2) |n,m,k,r,s,v\rangle &=& Q^{-(n-m+k+r-s+v)}|n,m,k,r,s-1,v-1\rangle
\nomb \\
 \pi(\Th_2^*) |n,m,k,r,s,v\rangle &=&
Q^{-(n-m+k+r-s+v)}|n,m,k,r,s+1,v+1\rangle
\eqa
\section{Conclusion}
We have given the algebraic embedding of two $q$-euclidean groups in three
dimensions
We have given an irreducible Hilbert space representation for the function
algebra of $I\SOd$.
With a reasoning along the lines of \cite{woro} it should be possible to prove
that $I\SOd$
exists on a $C^*$-algebra level.

\vspace{1cm}

\noi

{\bf Acknowledgments}

We would like to thank R. Engeldinger, J. Seifert and W. Weich for valuable
discussions.

\end{document}